\documentclass{conm-p-l}
\usepackage[all]{xy}

\theoremstyle{plain}
\newtheorem{theorem}{Theorem}[section]
\newtheorem{proposition}[theorem]{Proposition}
\newtheorem{lemma}[theorem]{Lemma}
\newtheorem{corollary}[theorem]{Corollary}
\newtheorem{conjecture}[theorem]{Conjecture}

\theoremstyle{definition}
\newtheorem{definition}[theorem]{Definition}
\newtheorem{example}[theorem]{Example}

\theoremstyle{remark}
\newtheorem{remark}[theorem]{Remark}

\newcommand{\secref}[1]{Section~\ref{#1}}
\newcommand{\thmref}[1]{Theorem~\ref{#1}}
\newcommand{\propref}[1]{Proposition~\ref{#1}}
\newcommand{\lemref}[1]{Lemma~\ref{#1}}
\newcommand{\corref}[1]{Corollary~\ref{#1}}
\newcommand{\conjref}[1]{Conjecture~\ref{#1}}

\newcommand{\exref}[1]{Example~\ref{#1}}

\newcommand{\defref}[1]{Definition~\ref{#1}}
\newcommand{\remref}[1]{Remark~\ref{#1}}

\begin{document}

\title{Subgroups of the Group of Self-Homotopy Equivalences}

\author{Martin Arkowitz}

\address{Department of Mathematics,
         Dartmouth College, Hanover NH 03755 U.~S.~A.}

\email{Martin.Arkowitz@Dartmouth.edu}

\author{Gregory Lupton}

\address{Department of Mathematics,
         Cleveland State University,
         Cleveland OH 44115
         U.~S.~A.}

\email{Lupton@math.csuohio.edu}

\author{Aniceto Murillo}

\address{Departmento de Algebra, Geometr\'\i a y Topolog\'\i a,
Universidad de M\'alaga,
Ap. 59, 29080 M\'alaga,
Spain}

\email{Aniceto@agt.cie.uma.es}

\date{\today}



\begin{abstract}
Denote by $\mathcal{E}(Y)$ the group of homotopy classes of
self-homotopy equivalences of a finite-dimensional complex $Y$. We
give a selection of results about certain subgroups of
$\mathcal{E}(Y)$. We establish a connection between the Gottlieb
groups of $Y$ and the subgroup of $\mathcal{E}(Y)$ consisting of
homotopy classes of self-homotopy equivalences that fix homotopy
groups through the dimension of $Y$, denoted by
$\mathcal{E}_{\#}(Y)$. We give an upper bound for the solvability
class of $\mathcal{E}_{\#}(Y)$ in terms of a cone decomposition of
$Y$.  We dualize the latter result to obtain an upper bound for
the solvability class of the subgroup of $\mathcal{E}(Y)$
consisting of homotopy classes of self-homotopy equivalences that
fix cohomology groups with various coefficients. We also show that
with integer coefficients, the latter group is nilpotent.
\end{abstract}

\maketitle

\section{Introduction and Preliminaries}\label{sec:intro}

Let $Y$ be a CW-complex of dimension $N$ and $\mathcal{E}(Y)$ the
group of homotopy classes of self-homotopy equivalences of $Y$. In
this paper we present a sample of results about a number of
subgroups of $\mathcal{E}(Y)$.  We denote by $\mathcal{E}_{\#}(Y)$
the following proto-typical such subgroup:
\begin{displaymath}
\mathcal{E}_{\#}(Y) = \{ f \in \mathcal{E}(Y) \mid f_{\#} = 1
\colon \pi_i(Y) \to \pi_i(Y),\   \text{for all $i\leq N$} \}.
\end{displaymath}
In \secref{sec:Gottlieb}, we give a way to construct elements in
$\mathcal{E}_{\#}(Y)$.  This is of interest since it provides a
connection between the \emph{Gottlieb group} of $Y$ and certain
subgroups of self-homotopy equivalences
(\thmref{thm:IntegralTheta}). Next, in \secref{sec:solvEsharp}, we
consider questions about the solvability and nilpotency of
$\mathcal{E}_{\#}(Y)$. For example, we show that if $Y$ is the
cofibre of a map between two wedges of spheres, then
$\mathcal{E}_{\#}(Y)$ is an abelian group
(\corref{cor:two-stage}). This result generalizes into a simple
upper bound on the solvability of $\mathcal{E}_{\#}(Y)$ in terms
of a cone-length invariant of $Y$ (\thmref{thrm:solvable}). In
\secref{sec:Estar} we dualize these results to obtain upper bounds
for the solvability of the group of equivalences that
 fix cohomology with different coefficients (\thmref{thrm:solvable-star}).  We also show
that the subgroup of self-homotopy equivalences which fix the
integral cohomology of a finite complex is  a nilpotent group
(\propref{prop:E* nilpotent}).

We now review some standard
material that we will use.  A cofibration sequence
\begin{displaymath}
\xymatrix{Z \ar[r]^{\gamma} & Y \ar[r]^j & X \ar[r]^-{q} &
\Sigma Z ,\\}
\end{displaymath}
where $X$ is the mapping cone of $\gamma$, gives a homotopy
coaction $c\colon X\to X\vee\Sigma Z$,  obtained by pinching the
`equator' of the cone of $Z$ to a point. The coaction induces an
action of $[\Sigma Z,W]$ on $[X,W]$ for any space $W$
(cf.~\cite[Chap.\,15]{Hil65} for details).  This is defined as
follows: If $\alpha\in[\Sigma Z,W]$ and $f\in[X,W]$, then
$f^\alpha$ is the composition
\begin{displaymath}
\xymatrix{X \ar[r]^-{c} & X \vee \Sigma Z\ar[r]^{f \vee \alpha} & W\vee
W \ar[r]^-{\nabla}& W.\\}
\end{displaymath}
The following properties of this action are well-known, and follow
easily from the definitions:
\begin{enumerate}
\item If $h \colon W \to W'$, then $h(f^{\alpha}) =
(hf)^{h\alpha}$. \item If $\alpha, \beta \in [\Sigma Z,W]$, then
$(f^{\alpha})^{\beta} = f^{(\alpha + \beta)}$.
\end{enumerate}
Next, consider the following portion of the Puppe sequence
associated to the above cofibration sequence:
\begin{displaymath}
\xymatrix{ [\Sigma Z, W] \ar[r]^{q^*}& [X, W]
\ar[r]^{j^*}& [Y, W]. \\ }
\end{displaymath}
As is also well-known, the orbits of the action are precisely the
pre-images of $j^*$.  That is, for $f, g \in [X, W]$, we have
$f^\alpha=g$ for some $\alpha\in[\Sigma Z,W]$ if and only if
$fj=gj$.

Next, we review some notation and terminology for groups.  Suppose that $G$
is a group and $H$ and $K$ are subgroups.  Then $H \triangleleft G$ denotes that
$H$ is normal in $G$ and $[H, K]$ denotes the subgroup generated by
commutators of elements of $H$ with elements of $K$. A
\emph{normal chain} for $G$ is a sequence of subgroups
\begin{displaymath}
G = G_1 \supseteq G_2 \supseteq \cdots \supseteq G_{k+1} \supseteq
\cdots
\end{displaymath}
with $G_{i+1} \triangleleft G_{i}$ for $i \geq 1$. If $[G_{i},
G_{i}] \subseteq G_{i+1}$ for each $i$, then the sequence is
called a \emph{solvability series}.  If, further, $G_{k+1} =
\{1\}$, then we say that $G$ is \emph{solvable of class} $\leq k$
and write $\mathrm{solv}\, G \leq k$. Analogously, given a normal
chain as above with each $G_{i} \triangleleft G$ and $[G, G_{i}]
\subseteq G_{i+1}$, then it is called a \emph{nilpotency series}.
In this case, if $G_{k+1} = \{1\}$, then we say that $G$ is
\emph{nilpotent of class} $\leq k$ and write $\mathrm{nil}\, G
\leq k$.  Clearly, we have $\mathrm{nil}\, G \geq \mathrm{solv}\,
G$.  In addition,  we write the identity homomorphism of a group
and the trivial homomorphism between two groups as $1 \colon G \to
G$ and $0 \colon G \to H$, respectively. This notation is also
used for sets with a distinguished element.

Finally, we fix our topological conventions
and notation. By a \emph{space}, we mean a connected CW-complex
of finite type.
Usually, we will be interested in
finite-dimensional CW-complexes.  When we discuss rational spaces,
we will specialize to $1$-connected CW-complexes.
As is well known, such a space $X$ admits a
\emph{rationalization}, which is denoted by $X_{\mathbb{Q}}$.
Similarly, a map of $1$-connected finite complexes $f \colon X \to
Y$ admits a rationalization map $f_{\mathbb{Q}} \colon X_{\mathbb{Q}}
\to Y_{\mathbb{Q}}$. A general reference for rationalization is
\cite{H-M-R}.  Furthermore,
we do
not distinguish notationally between a map and its homotopy class.
We
write $X \equiv Y$ to denote that the spaces $X$
and $Y$ have the same homotopy type. The
identity map of a space $X$ is denoted
$\iota \colon X \to X$ and the trivial map between two spaces $*
\colon X \to Y$.

\section{A Connection with the Gottlieb Group}\label{sec:Gottlieb}

We consider the situation as in \secref{sec:intro} of a mapping
cone sequence
\begin{displaymath}
\xymatrix{Z \ar[r]^{\gamma} & Y \ar[r]^j & X \ar[r]^-{q} &
\Sigma Z ,\\}
\end{displaymath}
and the induced action
of $[\Sigma Z,W]$ on $[X,W]$
which yields $f^\alpha
\in [X,W]$
for $\alpha\in[\Sigma Z,W]$
and $f\in[X,W]$.  We are interested in the effect that $f^\alpha$
has on homology and homotopy groups.  This is
described in the following result.

\begin{proposition}\label{prop:iotaalpha*}
For the above cofibration sequence, suppose $f\in[X,W]$ and
$\alpha\in[\Sigma Z,W]$.
\begin{enumerate}
\item The induced homology homomorphism $(f^{\alpha})_* \colon
H_i(X) \to H_i(W)$ is given by $(f^{\alpha})_*(x) = f_*(x) +
\alpha_* q_*(x)$, for each $x \in H_i(X)$.
\item Suppose that $(f,\alpha)\colon X\vee\Sigma Z\to W$ factors
through the product $X\times\Sigma Z$.  Then the induced homotopy
homomorphism $(f^{\alpha})_\# \colon \pi_i(X) \to \pi_i(W)$ is
given by $(f^{\alpha})_\#(x) = f_\# (x) + \alpha_\# q_\#(x)$, for
each $x \in \pi_i(X)$.
\end{enumerate}
\end{proposition}

\begin{proof}
(1) This follows directly from the commutative diagram
\begin{displaymath} \xymatrix{H_i(X) \ar[r]^-{c_*} \ar[rd]_-{(1, q_*)}&
H_i(X \vee \Sigma Z)
 \ar[r]^-{(f \vee\alpha)_*}\ar[d]^{(p_{1*}, p_{2*})} & H_i(W \vee W)
 \ar[r]^-{\nabla_*}\ar[d]_{(p_{1*}, p_{2*})} & H_i(W)\\
   & H_i(X) \oplus H_i(\Sigma Z) \ar[r]_-{f_*\oplus\alpha_*} & H_i(W) \oplus
   H_i(W)\ar[ru]_{+} & \\}
\end{displaymath}
in which the vertical
maps are isomorphisms induced by the two projections $p_1$ and $p_2$
and the top row is the
homomorphism induced by $f^{\alpha}$.

(2) Let $\sigma \colon S^i \to S^i \vee S^i$ denote the standard
comultiplication. Write $f_\#(x) + \alpha_\# q_\#(x)$ as the
composition
\begin{displaymath} \xymatrix{S^{i}
\ar[r]^-{(x\vee qx)\sigma}& X \vee \Sigma Z \ar[r]^-{(f,\alpha)} & W\\}
\end{displaymath}
and $(f^{\alpha})_\#(x)$ as the composition
\begin{displaymath} \xymatrix{S^{i}
\ar[r]^-{cx}& X \vee \Sigma Z \ar[r]^-{(f,\alpha)} & W\\}.
\end{displaymath}
By hypothesis, we can factor $(f,\alpha)$ through the product as
$(f,\alpha) = a\circ j \colon X \vee \Sigma Z \to X \times \Sigma
Z \to W$, for some $a \colon X \times \Sigma Z \to W$.  It is
straightforward to prove that $j(x\vee qx)\sigma = jcx \colon S^i
\to X \times \Sigma Z$, by checking that their projections onto
each summand are homotopic.
\end{proof}

We now specialize to a mapping cone sequence of the form
\begin{displaymath}
\xymatrix{S^{n-1} \ar[r]^{\gamma} & Y \ar[r]^j & X \ar[r]^-{q} &
\Sigma S^{n-1} \equiv S^{n} ,\\}
\end{displaymath}
i.e., $X = Y\cup_{\gamma}e^n$.  Then we have an action of
$\pi_n(X)$ on $[X, X]$.  We consider elements of the form
$\iota^\alpha \in [X,X]$, where $\iota$ is the identity map of $X$
and $\alpha \in \pi_n(X)$.  In general, these maps are not
self-homotopy equivalences.  However, by adding certain
hypotheses, we obtain maps in $\mathcal{E}_{\#}(X)$, or some other
subgroup of $\mathcal{E}(X)$. This approach is similar to that
taken in \cite{Au-Le}, but instead of assuming $j_{\#} \colon
\pi_{\#}(Y) \to \pi_{\#}(X)$ is onto, we shall consider
restrictions on the homotopy element $\alpha \in \pi_n(X)$. Recall
that the $n$'th \emph{Gottlieb group} of $X$, denoted $G_n(X)$,
consists of those $\alpha \in \pi_n(X)$ for which there is an
associated map $a \colon X \times S^n \to X$ such that the
following diagram commutes:
\begin{displaymath} \xymatrix{X \vee S^n
\ar[r]^-{(\iota ,\alpha)}\ar[d]_{j} & X
\\ X \times S^n \ar[ru]_{a}\\}.
\end{displaymath}
See \cite{Gott} for various results on the
groups $G_n(X)$.

Next, we introduce another subgroup of $\mathcal{E}(X)$. Define
\begin{displaymath}
\mathcal{E}_*(X) = \{ f \in \mathcal{E}(X) \mid f_* = 1 \colon
H_i(X) \to H_i(X), \  \text{for all $i$} \}.
\end{displaymath}
We apply \propref{prop:iotaalpha*} and obtain the following
consequence.

\begin{corollary}\label{cor:alpha*-q*}
Let $X = Y\cup_{\gamma}e^n$ be a $1$-connected CW-complex and
$\alpha \in \pi_n(X)$.
\begin{enumerate}
\item   $\iota^{\alpha} \in \mathcal{E}_{*}(X)$ if
and only if $\alpha_* q_* = 0 \colon H_{n}(X) \to H_{n}(X)$.

\item Suppose $\alpha \in G_n(X)$, and $X$ is of dimension $N$.
Then $\iota^{\alpha} \in \mathcal{E}_{\#}(X)$ if and only if
$\alpha_\# q_\# = 0 \colon \pi_{i}(X) \to \pi_{i}(X)$ for $i \leq
N$.
\end{enumerate}
\end{corollary}

\begin{proof}
It is immediate from \propref{prop:iotaalpha*} that $\iota^\alpha$
induces the identity on homology groups.  Since $X$
is a 1-connected CW-complex, $\iota^\alpha$ is a homotopy equivalence.  Hence
$\iota^{\alpha} \in \mathcal{E}_{*}(X)$.  This establishes (1),
and (2) follows similarly.
\end{proof}

Hence, we are interested in finding situations in which $\alpha q
\colon X \to X$ induces the trivial homomorphism, either in
homology or homotopy. Our first result is an
integral result. Following this, we shall focus on the rational
setting, where
more information can be obtained.

\begin{theorem}\label{thm:IntegralTheta}
Let $X = Y\cup_{\gamma}e^n$ be a $1$-connected $n$-dimensional
complex. Suppose that $q_\# = 0 \colon \pi_{n}(X) \to
\pi_{n}(S^n)$.  Then there is a homomorphism
\begin{displaymath}
\Theta \colon G_n(X) \to \mathcal{E}_{\#}(X),
\end{displaymath}
defined by $\Theta(\alpha) = \iota^{\alpha}$ for $\alpha \in
G_n(X)$.  This homomorphism restricts to
\begin{displaymath}
\Theta' \colon G_n(X)\cap \mathrm{ker}\,h_n \to
\mathcal{E}_{*}(X)\cap\mathcal{E}_{\#}(X),
\end{displaymath}
where $h_n \colon \pi_n(X) \to H_n(X)$ denotes the Hurewicz
homomorphism.
\end{theorem}

\begin{proof}
Let $\alpha \in G_n(X)$.  Since $\pi_{i}(S^n)=0$ for $i<n$, the
hypothesis gives that $q_{\#} = 0 \colon \pi_{i}(X) \to
\pi_{i}(S^n)$ for $i \leq n$. Hence, by \corref{cor:alpha*-q*},
$\iota^{\alpha} \in \mathcal{E}_{\#}(X)$. Now suppose $\beta$ is
any element in $\pi_n(X)$.  Since $\iota^{\alpha} \in
\mathcal{E}_{\#}(X)$ and $\beta \in \pi_n(X)$, we have that
$\iota^{\alpha}(\beta)=\beta$.  Thus, by the properties of the
action listed in the introduction, we have
\begin{displaymath}
\iota^{\alpha}\iota^{\beta} =
(\iota^{\alpha})^{\iota^{\alpha}(\beta)} = \iota^{\alpha +
\iota^{\alpha}(\beta)}=\iota^{\alpha+\beta}.
\end{displaymath}
Therefore $\Theta$ is a homomorphism.

Now suppose $\alpha$ is any element in $\mathrm{ker}\,h_n$.  Then
 $\alpha_* \colon H_n(S^n) \to
H_n(X)$ is zero. Since $H_i(S^n) = 0$ for positive $i \neq n$,
\corref{cor:alpha*-q*} implies that $\iota^{\alpha} \in
\mathcal{E}_{*}(X)$.  Thus $\Theta$ restricts to $\Theta'$ as
claimed.
\end{proof}

\begin{remark}
It is known that $G_n(X) \subseteq \mathrm{ker}\,h_n$ under
certain hypotheses (cf. \cite[Th.4.1]{Gott}), so the homomorphism
$\Theta$ and its restriction $\Theta'$ may agree.
\end{remark}

We illustrate \thmref{thm:IntegralTheta} with an example.

\begin{example}\label{ex:S2S3}
Take $X = S^2 \times S^3=S^2\vee S^3\cup_{[i_1,i_2]} e^5$.  This
kind of example has been considered previously (cf. \cite{Ark-Mar,
Saw}), but here we put it into the context discussed
above.

As is well-known, $S^3$ is an $H$-space and therefore satisfies
$G_i(S^3) = \pi_i(S^3)$ for all $i$.  Further, the Gottlieb group
preserves products so $G_5(X) = G_5(S^2) \oplus G_5(S^3)$.  Since
$\pi_5(S^3) = \mathbb{Z}_2$, there is at least a non-trivial
element of order $2$ in $G_5(X)$.  Next, consider $q_{\#} \colon
\pi_{5}(X) \to \pi_{5}(S^5)$.  Since $\pi_5(X) = \pi_5(S^3) \oplus
\pi_5(S^2)$ is a finite group and $\pi_{5}(S^5)$ is infinite
cyclic, it follows that $q_{\#}$ is zero in this dimension. From
\thmref{thm:IntegralTheta}, $\Theta$ defines a homomorphism from
$G_5(X)$ to $\mathcal{E}_{\#}(X)$.  Notice that $h_5 \colon
\pi_5(X) \to H_5(X)$ is zero, since $H_5(X)$ is infinite cyclic.
 Therefore, $G_5(X)
\subseteq \mathrm{ker}\,h_5$ and $\Theta = \Theta' \colon G_5(X)
\to \mathcal{E}_{*}(X)\cap\mathcal{E}_{\#}(X)$.
\end{example}

Note that the homomorphism $\Theta$ in \thmref{thm:IntegralTheta}
may have trivial image in $\mathcal{E}_{\#}(X)$. For
\exref{ex:S2S3}, it follows from the computations in \cite[\S
6]{Ark-Mar} that $\Theta$ is actually injective.  However, it
seems to be difficult to give general conditions to guarantee that
$\Theta$ is injective.  Rather than doing this by placing strong
hypotheses on our spaces, we turn now to the rational setting. For
a 1-connected CW-complex $X$, the \emph{rational Gottlieb group}
of $X$ is the Gottlieb group of the rationalization of $X$, that
is, $G_n(X_{\mathbb{Q}})$. Notice that by \cite{Lang}, we have
$G_n(X_{\mathbb{Q}}) \cong G_n(X)\otimes \mathbb{Q}$ for each $n$.
In contrast to the ordinary Gottlieb groups, much is known about
the rational Gottlieb groups by results of F{\'e}lix-Halperin
\cite{Fe-Ha82}.  For instance, a $1$-connected, finite complex has
no non-trivial rational Gottlieb groups of even degree, and has
only finitely many non-trivial rational Gottlieb groups of odd
degree (see \cite{Fel89} for details).  We only touch on these
ideas here and avoid heavy use of rational techniques.

\begin{lemma}\label{lem:jsharpsurj}
Suppose we have a mapping cone sequence
\begin{displaymath}
\xymatrix{S^{n-1} \ar[r]^{\gamma} & Y \ar[r]^j & X \\}
\end{displaymath}
in which $Y$ and $X$ are $1$-connected.  If $\gamma_{\mathbb{Q}}
\neq
* \colon (S^{n-1})_{\mathbb{Q}} \to Y_{\mathbb{Q}}$, then
$(j_{\mathbb{Q}})_{\#} \colon \pi_{n}(Y_{\mathbb{Q}})$ $\to
\pi_{n}(X_{\mathbb{Q}})$ is surjective.
\end{lemma}

\begin{proof}
This can be argued using the long exact homotopy and homology
sequences, together with the relative Hurewicz theorem, for the
pair $(X, Y)$. Alternatively, Quillen minimal models can be used.
We omit the details.
\end{proof}

\begin{remark}
Notice we assert that $(j_{\mathbb{Q}})_{\#}$ is onto in degree
$n$ only, and not in all degrees.  In the latter case, the cell
attachment is called an \emph{inert} cell attachment \cite{Ha-Le}.
This is one of the hypotheses used in \cite{Au-Le}, but it is not
satisfied by some of the examples we have in mind.
\end{remark}

We will see that under our hypotheses, rational equivalences of the form
$\iota^{\alpha}$ are contained in a smaller subgroup of
$\mathcal{E}(X)$ than $\mathcal{E}_{\#}(X)$. We introduce the following
notation: For $r \leq \infty$, define
\begin{displaymath}
\mathcal{E}_{\#r}(X) = \{ f \in \mathcal{E}(X) \mid f_{\#} = 1
\colon \pi_i(X) \to \pi_i(X),\ \text{for all $i \leq r$} \}.
\end{displaymath}
Note that $f\in \mathcal{E}_{\# \infty}(X)$ if and only if
$f$ induces the identity homomorphism of \emph{all}
homotopy groups.

The following is our basic rational result.

\begin{theorem}\label{thm:RationalPhi}
Let $X = Y\cup_{\gamma}e^n$ be a $1$-connected CW complex with $n$
odd and $\gamma_{\mathbb{Q}} \neq * \colon (S^{n-1})_{\mathbb{Q}}
\to Y_{\mathbb{Q}}$. Then there is a homomorphism
\begin{displaymath}
\Phi \colon G_n(X_{\mathbb{Q}}) \to
\mathcal{E}_{\#\infty}(X_{\mathbb{Q}}),
\end{displaymath}
defined by $\Phi(\alpha) = \iota^{\alpha}$ for $\alpha \in
G_n(X_{\mathbb{Q}})$.  This homomorphism restricts to
\begin{displaymath}
\Phi' \colon G_n(X_{\mathbb{Q}}) \cap \mathrm{ker}\,h_n \to
\mathcal{E}_{*}(X_{\mathbb{Q}}) \cap
\mathcal{E}_{\#\infty}(X_{\mathbb{Q}}),
\end{displaymath}
where $h_n$ denotes the rational Hurewicz homomorphism $h_n \colon
\pi_n(X_{\mathbb{Q}}) \to H_n(X_{\mathbb{Q}})$.
\end{theorem}

\begin{proof}
We proceed as in the proof of \thmref{thm:IntegralTheta}.  First,
we claim that $(q_{\mathbb{Q}})_{\#} = 0 \colon
\pi_i(X_{\mathbb{Q}}) \to \pi_i(S^n_{\mathbb{Q}})$, for all $i$.
Since $n$ is odd, we have
$\pi_i(S^n_{\mathbb{Q}}) = 0$ for $i\ne n$.
Hence we must only check that $(q_{\mathbb{Q}})_{\#} = 0$ in
degree $n$.  By \lemref{lem:jsharpsurj},
$(j_{\mathbb{Q}})_{\#} \colon \pi_n(Y_{\mathbb{Q}}) \to
\pi_n(X_{\mathbb{Q}})$ is surjective. Given $x \in
\pi_n(X_{\mathbb{Q}})$, write $x = (j_{\mathbb{Q}})_{\#}(y)$, for some
$y \in \pi_n(Y_{\mathbb{Q}})$.  Then $(q_{\mathbb{Q}})_{\#}(x) =
(q_{\mathbb{Q}})_{\#}(j_{\mathbb{Q}})_{\#}(y) =
\big((qj)_{\mathbb{Q}}\big)_{\#}(y)=0$ since $qj = * $, and the claim
follows.

Now a simple modification of the proof of
\lemref{cor:alpha*-q*} yields that $\iota^\alpha \in
\mathcal{E}_{\#\infty}(X_{\mathbb{Q}})$, for each $\alpha \in
G_n(X_{\mathbb{Q}})$.  The remainder of the argument follows
exactly as in the proof of \thmref{thm:IntegralTheta}.
\end{proof}

\begin{remark}
Notice that, unlike \thmref{thm:IntegralTheta}, there is no
restriction on the dimension of $X$ in \thmref{thm:RationalPhi},
and the attached cell need not be top-dimensional.  If $X$ is a
$1$-connected, finite complex, then there is no generality lost in
assuming $n$ odd since a $1$-connected, finite complex has
no non-trivial rational Gottlieb groups of even degree.
\end{remark}

Although \thmref{thm:RationalPhi} is a rational result, we are
able to `de-rationalize' it to obtain the following integral
consequence.

\begin{theorem}
Let $X = Y\cup_{\gamma}e^n$ be a $1$-connected finite complex with
$n$ odd and $\gamma \in \pi _{n-1}(Y)$ not of finite order. If the
homomorphism $\Phi$ from \thmref{thm:RationalPhi} is non-zero,
then for each $r$ with $\mathrm{dim}\,X \leq r < \infty$, there
are elements of infinite order in $\mathcal{E}_{\#r}(X)$.
\end{theorem}

\begin{proof}
Suppose $\Phi(\alpha) = \iota ^{\alpha}$ is not the identity element in
$\mathcal{E}_{\#\infty}(X_{\mathbb{Q}})$.  Since this latter is a
$\mathbb{Q}$-local group, it contains no non-trivial elements of finite order,
and hence $(\iota^\alpha)^k \not= \iota \colon X_{\mathbb{Q}} \to
X_{\mathbb{Q}}$, for all $k$.  By \cite{Mar}, we have
$\mathcal{E}_{\#r}(X_{\mathbb{Q}}) \cong
\big(\mathcal{E}_{\#r}(X)\big) _{\mathbb{Q}}$, for each $r$ with
$\text{dim}\,X \leq r < \infty$. From this it follows that for
each $r$ there is some positive integer $p$ and some element $f
\in \mathcal{E}_{\#r}(X)$ such that $f_{\mathbb{Q}} =
(\iota^\alpha)^p$.  Since $f_{\mathbb{Q}}$ is of infinite order in
$\mathcal{E}_{\#r}(X_{\mathbb{Q}})$, the same must be true of $f$
in $\mathcal{E}_{\#r}(X)$.
\end{proof}

\section{Solvability of $\mathcal{E}_{\#}(Y)$}\label{sec:solvEsharp}

A result of Dror-Zabrodsky asserts that if $Y$ is a finite
complex, then $\mathcal{E}_{\#}(Y)$ is a nilpotent group
\cite{D-Z}. One can ask, therefore, whether there are reasonable
estimates for the nilpotency, or perhaps the solvability, of
$\mathcal{E}_{\#}(Y)$ in terms of the usual algebraic topological
invariants of $Y$.   Several results have been established that
relate the nilpotency or solvability of $\mathcal{E}_{\#}(Y)$, or
some similar group, to the \emph{Lusternik-Schnirelmann category}
of $Y$, or related invariants (cf. \cite{Ark-Lup, Fe-Mu97,
Fe-Mu98, Sc-Ta99}). Some of these apply in a rational setting, and
others in an integral setting. Typically, these results give an
upper bound on the nilpotency or solvability of the group.

We begin by discussing a topological invariant which appears in our results.

\begin{definition}\label{def:scl}
For any space $X$, a \emph{spherical cone decomposition of} $X$
\emph{of length} $n$, is a sequence of cofibrations
\begin{displaymath}
\xymatrix{L_i \ar[r]^{\gamma_i} & X_i \ar[r]^{j_i} & X_{i+1},\\}
\end{displaymath}
for $0 \leq i < n$, such that each $L_i$ is a finite wedge of
spheres, $X_0$ is contractible and $X_n \equiv X$.  We define the
\emph{spherical cone-length} of $X$, denoted by $\mathrm{scl}(X)$,
as follows:  If $X$ is contractible, then set $\mathrm{scl}(X) =
0$. Otherwise, $\mathrm{scl}(X)$ is the smallest positive integer
$n$ such that there exists a spherical cone decomposition of $X$
of length $n$. If no such integer exists, set $\mathrm{scl}(X) =
\infty$. If $X$ is a finite-dimensional complex and we have a
spherical cone decomposition of $X$ of length $n$ in which, in
addition, $\text{dim}\,L_i < \text{dim}\,X$ for $i = 0, \ldots,
n-1$, then this is called a \emph{restricted spherical cone
decomposition of} $X$ \emph{of length} $n$.  We then define the
\emph{restricted spherical cone-length} of $X$, denoted
$\mathrm{rscl}(X)$, using only restricted spherical cone
decompositions in place of ordinary spherical cone decompositions.
\end{definition}

\begin{remark}\label{rem:scl}
Spherical cone-length has been considered in \cite{Co94, Sc-Ta99}.
If we denote the Lusternik-Schnirelmann category of $X$ by
$\mathrm{cat}(X)$, then it is known that $\mathrm{cat}(X) \leq
\mathrm{scl}(X)$.  Note that a space $X$ with $\mathrm{scl}(X) =
1$ is homotopy equivalent to a wedge of spheres and that a space
$X$ with $\mathrm{scl}(X) \leq 2$ is homotopy equivalent to the
cofibre of a map between wedges of spheres. Furthermore, the
cell-structure of a finite-dimensional complex $X$ provides a
restricted spherical cone decomposition of length $\leq$ the
number of dimensions in which there are positive-dimensional
cells.
\end{remark}

We introduce a bit more notation before proving the main result of this
section. Once again, $Y$ is a complex of dimension
$N$.  We define
\begin{align*}
\mathcal{E}_{k}(Y) = \{ f \in \mathcal{E}(Y) \mid f_* = 1 \colon
[X, Y] \to &[X, Y],\   \text{for every complex $X$} \\ &
\text{with dim $X \leq N$ and $\mathrm{rscl}(X) \leq k$}\}.
\end{align*}
In particular, we have $\mathcal{E}_{1}(Y) = \mathcal{E}_{\#}(Y)$.
Also, there is a chain of subgroups
\begin{equation}\label{eqn:sequence}
\mathcal{E}_{\#}(Y) = \mathcal{E}_{1}(Y) \supseteq
\mathcal{E}_{2}(Y) \supseteq \cdots \supseteq \mathcal{E}_{k}(Y)
\supseteq \cdots.
\end{equation}
Clearly we have $\mathcal{E}_{k}(Y) \triangleleft
\mathcal{E}_{k-1}(Y)$: For if $f \in \mathcal{E}_{k}(Y)$, $g \in
\mathcal{E}_{k-1}(Y)$, $\mathrm{dim}\,X \leq N$ and
$\mathrm{rscl}(X) \leq k$, then $f_*g^{-1}_* = g^{-1}_* \colon [X,
Y] \to [X, Y]$.  Hence
\begin{displaymath}
(gfg^{-1})_* = g_* f_* g^{-1}_*  = g_* g^{-1}_* = 1,
\end{displaymath}
and so $gfg^{-1} \in \mathcal{E}_{k}(Y)$.  Therefore, the series
(\ref{eqn:sequence}) is a normal chain. Furthermore, if
$\mathrm{rscl}(Y) \leq k$, then $\mathcal{E}_{k}(Y) = 1$.  Then
 we have a normal chain
\begin{equation*}
\mathcal{E}_{\#}(Y)  \supseteq \mathcal{E}_{2}(Y) \supseteq \cdots
\supseteq \mathcal{E}_{k}(Y) = \{1\}.
\end{equation*}

The following is the main result of this section.

\begin{theorem}\label{thrm:solvable}
The series (\ref{eqn:sequence}) is a solvability series, i.e.,
$[\mathcal{E}_{i}(Y), \mathcal{E}_{i}(Y)] \subseteq
\mathcal{E}_{i+1}(Y)$.  Consequently, we have
\begin{displaymath}
\mathrm{solv}\, \mathcal{E}_{\#}(Y)\leq \mathrm{rscl}(Y) - 1.
\end{displaymath}
\end{theorem}

\begin{proof}
Let $f, g \in \mathcal{E}_{i}(Y)$ and $X$ be a complex with
$\mathrm{rscl}(X) = i+1$ and dim$\,X \le N$.  It suffices to show
that $f_*g_*(h) = g_*f_*(h)$ for every $h \in [X, Y]$.  Consider
the last cofibre sequence in a length-$(i+1)$ restricted spherical
cone decomposition of $X$,
\begin{displaymath}
\xymatrix{L_i \ar[r]^{\gamma} & X_i \ar[r]^{j} & X_{i+1}, \\}
\end{displaymath}
where $L_i$ is a wedge of spheres and $X_{i+1} \equiv X$.  Now,
since $f \in \mathcal{E}_{i}(Y)$, it follows that $j^*(fh) =
f_*(hj) = j^*(h)$. Thus, from the properties of the coaction
reviewed in \secref{sec:intro}, there is some $\alpha \in [\Sigma
L_i, Y]$ such that $fh = h^{\alpha}$. Similarly, there is some
$\beta \in [\Sigma L_i, Y]$ such that $g_*(h) = h^{\beta}$. Note
also that $f\beta = \beta$ since $f \in \mathcal{E}_{i}(Y)
\subseteq \mathcal{E}_{1}(Y)$, and $\Sigma L_i$ is a wedge of
spheres of dimension $\leq N$. Now we have
\begin{displaymath}
f_*g_*(h) = f(h^{\beta}) = (fh)^{f\beta} = (h^{\alpha})^{\beta} =
h^{\alpha+\beta}.
\end{displaymath}
A similar computation yields $g_*f_*(h) = h^{\beta + \alpha}$.
Since $[\Sigma L_i, Y]$ is abelian, the proof is
complete.\end{proof}

\begin{remark}
A result analogous to \thmref{thrm:solvable} has been proved by
Scheerer and Tanr{\'e} in \cite[Th.6]{Sc-Ta99}.  We note the
differences and similarities between these results. Theorem 6 in
\cite{Sc-Ta99} is proved for the group of equivalences of a space
$Y$ relative to certain fixed classes of spaces (though our proof
could be easily modified to hold for these classes).  When the
class consists of wedges of spheres, the corresponding group of
equivalences is $\mathcal{E}_{\#\infty}(Y)$. The upper bound for
the solvability of the group of equivalences
$\mathcal{E}_{\#\infty}(Y)$ relative to the class of wedges of
spheres given in \cite{Sc-Ta99} is then the so-called spherical
category of $Y$, which is less than or equal to the spherical
cone-length of $Y$ minus one. On the other hand, the group of
equivalences $\mathcal{E}_{\#}(Y)$ that we consider in
\thmref{thrm:solvable} is larger than $\mathcal{E}_{\#\infty}(Y)$.
Furthermore, the two proofs are similar, but the solvability
series in \thmref{thrm:solvable} appears to be different from the
one in \cite[Th.6]{Sc-Ta99}.
\end{remark}

\thmref{thrm:solvable} easily gives the next two corollaries.

\begin{corollary}\label{cor:two-stage}
If $\mathrm{rscl}(Y) \leq 2$, that is, $Y$ is the cofibre of a map
between wedges of spheres, then $\mathcal{E}_{\#}(Y)$ is abelian.
\end{corollary}

\begin{corollary}\label{cor:number of cells}
For $Y$ any finite-dimensional complex, $\mathcal{E}_{\#}(Y)$ is
solvable, with $\mathrm{solv}\, \mathcal{E}_{\#}(Y)\leq k-1$,
where $k$ is the number of dimensions in which there are
positive-dimensional cells.
\end{corollary}

\begin{proof}
This follows by \remref{rem:scl}.
\end{proof}

We can modify much of the previous material to deal with
equivalences that fix all homotopy groups.  This also allows us to
deal with the case in which $Y$ is an arbitrary space which is not
necessarily a finite-dimensional complex.  Define
\begin{displaymath}
\mathcal{E}'_{k}(Y) = \{ f \in \mathcal{E}(Y) \mid f_* = 1 \colon
[X, Y] \to [X, Y], \ \text{for all $X$ with}\ \mathrm{scl}(X) \leq
k\}.
\end{displaymath}
Then there is a normal chain
\begin{equation}\label{eqn:sequence-dashed}
\mathcal{E}_{\#\infty}(Y) = \mathcal{E}'_{1}(Y) \supseteq
\mathcal{E}'_{2}(Y) \supseteq \cdots \supseteq\mathcal{E}'_{k}(Y)
\supseteq \cdots.
\end{equation}
Now the proof of \thmref{thrm:solvable} yields the following
analogous results.

\begin{theorem}{\cite[Th.\,6]{Sc-Ta99}}\label{thrm:solvable-dashed}
The series (\ref{eqn:sequence-dashed}) is a solvability series.
Therefore,
\begin{displaymath}
\mathrm{solv}\, \mathcal{E}_{\#\infty}(Y)\leq \mathrm{scl}(Y) - 1.
\end{displaymath}
\end{theorem}

We note that from a spectral sequence of Didierjean
\cite{Didierjean}, one can also obtain a different upper bound on
the solvability of $\mathcal{E}_{\#\infty}(Y)$ in terms of the
cohomology of $Y$ with coefficients in the homotopy groups of $Y$.

\begin{corollary}\label{cor:two-stage-dashed}
If $\mathrm{scl}(Y) \leq 2$, then $\mathcal{E}_{\#\infty}(Y)$ is
abelian.
\end{corollary}
In view of the results in this section, together with the bounds
found in \cite{Ark-Lup, Fe-Mu97, Fe-Mu98, Sc-Ta99}, it is natural
to believe that, for a finite-dimensional complex $Y$, the
nilpotency class of $\mathcal{E}_{\#}(Y)$ is bounded above by
$\mathrm{scl}(Y) - 1$. We have not been able to prove this, and so
we leave it as a conjecture.

\begin{conjecture}\label{conj:central}
For a finite-dimensional complex $Y$,
\begin{displaymath}
\mathrm{nil}\, \mathcal{E}_{\#}(Y)\leq \mathrm{rscl}(Y) - 1 \qquad
\text{and} \qquad \mathrm{nil}\, \mathcal{E}_{\#\infty}(Y)\leq
\mathrm{scl}(Y) - 1.
\end{displaymath}
\end{conjecture}

We note that Scheerer and Tanr{\'e} have conjectured that
$\mathrm{nil}\, \mathcal{E}_{\# \infty}(Y)$ is bounded above by
the spherical category of $Y$ \cite[\S7, (6)]{Sc-Ta99}.

\conjref{conj:central} would be established by showing that each
of the series (\ref{eqn:sequence}) and (\ref{eqn:sequence-dashed})
is a nilpotency series.  A direct proof of this would also give an
independent proof of the Dror-Zabrodsky result on the nilpotency
of $\mathcal{E}_{\#}(Y)$.

\section{Equivalences that Fix Cohomology Groups}\label{sec:Estar}

In this section, we dualize some of the ideas of the previous
section.  Although we did not use homotopy groups with
coefficients there, we do use coefficients here, since this is
more common with cohomology.

Let $\mathcal{G}$ be a collection of abelian groups and $X$ be a
space.  Define
\begin{align*}
\mathcal{E}^*_{\mathcal{G}}(X) = \{ f \in \mathcal{E}(X) \mid f^*
= 1 \colon H^i(X; G) \to H^i(X; G), \  \text{for all $i$ and all
$G \in \mathcal{G}$} \}.
\end{align*}
The following cases are of special interest:
\begin{enumerate}
\item $\mathcal{G} = \{ \mathbb{Z} \}$.  We write
$\mathcal{E}^*_{\mathcal{G}}(X)$ as $\mathcal{E}^*(X)$ in this
case.

\item $\mathcal{G} =$ all cyclic groups.  Then $f \in
\mathcal{E}^*_{\mathcal{G}}(X)$ if and only if $f \in
\mathcal{E}(X)$ and $f^* = 1 \colon H^i(X; G) \to H^i(X; G)$ for
every finitely-generated abelian group $G$.  We write
$\mathcal{E}^*_{\mathcal{G}}(X)$ as
$\mathcal{E}^*_{\mathrm{fg}}(X)$ in this case.  Note that
$\mathcal{E}^*_{\mathrm{fg}}(X) \subseteq \mathcal{E}^*(X)$.
\end{enumerate}

Next we define a topological invariant that plays a role dual to
that of spherical cone-length in the previous section.

\begin{definition} For $\mathcal{G}$ a collection of abelian
groups, call an Eilenberg-Mac-Lane space $K(G, m)$, with $G \in
\mathcal{G}$, a \emph{$\mathcal{G}$-Eilenberg-MacLane space}. For
any space $X$, a \emph{$\mathcal{G}$-fibre decomposition of} $X$
\emph{of length} $n$, is a sequence of fibrations
\begin{displaymath}
\xymatrix{X_{i+1} \ar[r]^{j_i} & X_{i} \ar[r]^{p_i} & K_{i},\\}
\end{displaymath}
for $0 \leq i < n$, such that each $K_i$ is a finite product of
$\mathcal{G}$-Eilenberg-MacLane spaces, $X_0$ is contractible and
$X_n \equiv X$. We define the \emph{$\mathcal{G}$-fibre-length} of
$X$, denoted by $\mathcal{G}\text{-}\mathrm{fl}(X)$, by dualizing
\defref{def:scl} in a straightforward way.
\end{definition}

Note that a space $X$ with $\mathcal{G}\text{-}\mathrm{fl}(X) = 1$
is homotopy equivalent to a product of
$\mathcal{G}$-Eilenberg-MacLane spaces. A space $X$ with
$\mathcal{G}\text{-}\mathrm{fl}(X) \leq 2$ is homotopy equivalent
to the fibre of a map between products of
$\mathcal{G}$-Eilenberg-MacLane spaces. Note also that when we
mention a product of $\mathcal{G}$-Eilenberg-MacLane spaces, we
allow factors with homotopy groups in different dimensions, so
that a product of $\mathcal{G}$-Eilenberg-MacLane spaces is not
itself a $\mathcal{G}$-Eilenberg-MacLane space in general.

Now define subgroups of $\mathcal{E}(X)$ as follows:
\begin{displaymath}
\mathcal{E}^*_{\mathcal{G}, s}(X) = \{ f \in \mathcal{E}(X) \mid
f^* = 1 \colon [X, Y] \to [X, Y],\ \text{for all $Y$ with
$\mathcal{G}\text{-}\mathrm{fl} (Y) \leq s$} \}.
\end{displaymath}
Then there is a normal chain of subgroups
\begin{equation}\label{eqn:sequence-star}
\mathcal{E}^*_{\mathcal{G}}(X) = \mathcal{E}^*_{\mathcal{G}, 1}(X)
 \supseteq \mathcal{E}^*_{\mathcal{G}, 2}(X) \supseteq \cdots \supseteq
\mathcal{E}^*_{\mathcal{G}, s}(X) \supseteq \cdots.
\end{equation}
The proof of normality for (\ref{eqn:sequence-star}) is similar to
the proof of normality for (\ref{eqn:sequence}) above.

A straightforward dualization of the proof of
\thmref{thrm:solvable} yields the following result.

\begin{theorem}\label{thrm:solvable-star}
The series (\ref{eqn:sequence-star}) is a solvability series. Thus
\begin{displaymath}
\mathrm{solv}\, \mathcal{E}^*_{\mathcal{G}}(X) \leq
\mathcal{G}\text{-}\mathrm{fl}(X) - 1.
\end{displaymath}
\end{theorem}

In \corref{cor:number of cells}, we showed that the number of
dimensions in which $Y$ has cells may be used to estimate
the spherical cone length.  We now indicate briefly how the
preceding notions can be modified for the
dual result.

\begin{definition}
A space $X$ is called a \emph{Postnikov piece} if there is some
$N$ such that $\pi_i(X) = 0$ for all $i > N$.  The smallest such
$N$ is called the \emph{homotopical dimension} of $X$, and is
denoted $\mathrm{h}\text{-}\mathrm{dim}\,X$.
\end{definition}

For a Postnikov piece $X$, we define another subgroup of $\mathcal{E}(X)$ as
\begin{align*}
\mathcal{E}^{*\prime}_{\mathcal{G}}(X) = \{ f \in \mathcal{E}(X)
\mid f^* = 1 \colon H^i(X; G) \to H^i(X; G), \ &\text{for all $i
\leq \mathrm{h}\text{-}\mathrm{dim}\,X$} \\& \text{and all $G \in
\mathcal{G}$} \}.
\end{align*}
We then define a \emph{restricted $\mathcal{G}$-fibre
decomposition} of a Postnikov piece $X$, of length $n$, as above
but with the additional condition that
$\mathrm{h}\text{-}\mathrm{dim}\,K_i \leq
\mathrm{h}\text{-}\mathrm{dim}\,X + 1$, for all $i$. This yields
the \emph{restricted $\mathcal{G}$-fibre-length} of $X$, denoted
by $\mathrm{r}\text{-}\mathcal{G}\text{-}\mathrm{fl}(X)$.

\begin{remark}\label{rem:G-fibre length}
Let $X$ be a $1$-connected space  and $X^{(N)}$ be the $N$'th
Postnikov section of $X$.  Set $\mathcal{G} =$ all cyclic groups
and $s=$ the number of non-trivial homotopy groups of $X^{(N)}$.
Then by taking the Postnikov decomposition of $X^{(N)}$, we see
that
$\mathrm{r}\text{-}\mathcal{G}\text{-}\mathrm{fl}(X^{(N)})\leq s$.
\end{remark}

Now if $X$ is a Postnikov piece, define
\begin{align*}
\mathcal{E}^{*\prime}_{\mathcal{G},s}(X) = \{ f \in \mathcal{E}(X)
\mid & f^* = 1 \colon [X, Y] \to [X, Y], \ \text{for all Postnikov
pieces $Y$} \\ &\text{such that h-dim}\,Y \leq
\mathrm{h}\text{-}\mathrm{dim}\,X \ \text{and}\
\mathrm{r}\text{-}\mathcal{G}\text{-}\mathrm{fl} (Y)\leq s \}.
\end{align*}
Then once again we have a normal chain of subgroups
\begin{equation}\label{eqn:sequence-star-dashed}
\mathcal{E}^{*\prime}_{\mathcal{G}}(X) =
\mathcal{E}^{*\prime}_{\mathcal{G}, 1}(X)
 \supseteq \mathcal{E}^{*\prime}_{\mathcal{G}, 2}(X) \supseteq \cdots
\supseteq
\mathcal{E}^{*\prime}_{\mathcal{G}, s}(X) \supseteq \cdots.
\end{equation}
Here also, the proof of normality is similar to the proof of
normality for (\ref{eqn:sequence}) above. A further adaptation of
the proof of \thmref{thrm:solvable} gives the following result.
\begin{theorem}\label{thrm:solvable-star-dashed}
The series (\ref{eqn:sequence-star-dashed}) is a solvability
series. Thus ,
\begin{displaymath}
\mathrm{solv}\, \mathcal{E}^{*\prime}_{\mathcal{G}}(X) \leq
\mathrm{r}\text{-}\mathcal{G}\text{-}\mathrm{fl}(X) - 1,
\end{displaymath}
for a Postnikov piece $X$.
\end{theorem}

This leads to the dual of
\corref{cor:number of cells}.

\begin{corollary}\label{cor:solvE*fg}
If $X$ is a $1$-connected finite complex, then
$\mathcal{E}^{*}_{\mathrm{fg}}(X)$ is solvable.  In particular, if
$X$ has dimension $N$ and there are $s$ non-trivial homotopy
groups in degrees $\leq N$, then $\mathrm{solv}\,
\mathcal{E}^{*}_{\mathrm{fg}}(X) \leq s -1$.
\end{corollary}

\begin{proof}
Any map $f \colon X \to X$ induces a corresponding map $\theta(f)
\colon X^{(N)} \to X^{(N)}$ of $N$'th Postnikov sections.  This
gives us a homomorphism $\theta \colon
\mathcal{E}^*_{\mathcal{G}}(X) \to$
$\mathcal{E}^{*\prime}_{\mathcal{G}}(X^{(N)})$ which is one-one.
Therefore, $\mathrm{solv}\, \mathcal{E}^*_{\mathcal{G}}(X) \leq
\mathrm{solv}\, \mathcal{E}^{*\prime}_{\mathcal{G}}(X^{(N)})$. But
if $\mathcal{G}$ is the collection of all cyclic groups, then it
is a consequence of \remref{rem:G-fibre length} and
\thmref{thrm:solvable-star-dashed} that $\mathrm{solv}\,
\mathcal{E}^{*\prime}_{\mathcal{G}}(X^{(N)}) \leq s-1$.
\end{proof}

Next, we compare the subgroups $\mathcal{E}^*(X)$ and
$\mathcal{E}^*_{\mathrm{fg}}(X)$. Below, we give a simple example
to illustrate that these two subgroups are distinct in general.
First, however, we obtain conditions under which they agree.

\begin{proposition}
For any space $X$, if
\begin{displaymath}
\mathrm{Hom}\big( \mathrm{Tor}(H^{i+1}(X),G), H^i(X)\otimes G
\big) = 0
\end{displaymath}
for all $i$ and all finitely-generated groups $G$, then
$\mathcal{E}^*(X) = \mathcal{E}^*_{\mathrm{fg}}(X).$
\end{proposition}

\begin{proof}
If $f \in \mathcal{E}^*(X)$, then the universal coefficient
theorem (cf. \cite[Thm.\,5.10]{Spa66}) gives
 a commutative diagram
with exact rows
\begin{displaymath}
\xymatrix{0 \ar[r] & H^i(X)\otimes G \ar[r]^{\iota} \ar[d]^{1} &
H^i(X ;G) \ar[r]^-{\pi} \ar[d]^{\phi} &\mathrm{Tor}(H^{i+1}(X),G)
\ar[r]\ar[d]^{1}& 0 \\ 0 \ar[r] & H^i(X)\otimes G \ar[r]^{\iota} &
H^i(X ;G) \ar[r]^-{\pi} &\mathrm{Tor}(H^{i+1}(X),G) \ar[r]& 0\\}
\end{displaymath}
where the middle homomorphism $\phi$ can be either $f^*$ or the
identity 1.  Thus there is a homomorphism $\rho \colon
\mathrm{Tor}(H^{i+1}(X),G) \to H^i(X)\otimes G$ such that $f^* - 1
= \iota\rho\pi$.  By hypothesis, $\rho = 0$ and so $f^* = 1$.
Therefore, $f \in \mathcal{E}^*_{\mathrm{fg}}(X)$.
\end{proof}

The following example illustrates that $\mathcal{E}^*(X)$ and
$\mathcal{E}^*_{\mathrm{fg}}(X)$ may differ.

\begin{example}
Let $X$ be a Moore space $M(G, n)$, for $n \geq 2$ and $G$ any
infinite, finitely-generated abelian group with torsion. Then it
follows from results of \cite{Ark-Mar} that $\mathcal{E}^*(X) \neq
\mathcal{E}^*_{\mathrm{fg}}(X)$
\end{example}

In \corref{cor:solvE*fg}, we showed that for a $1$-connected,
finite-dimensional complex $X$,  $\mathcal{E}^*_{\mathrm{fg}}(X)$
is solvable.  This raises the question of whether or not the group
is nilpotent.  We conclude the paper by showing that
$\mathcal{E}^*(X)$, and therefore
$\mathcal{E}^*_{\mathrm{fg}}(X)$, is a nilpotent group.  We will
use the following notation: Suppose that a group $G$ acts on an
abelian group $A$.  We define inductively a decreasing sequence of
subgroups of $A$ by setting $\Gamma^G_1(A)=A$ and $\Gamma^G_i(A)$
is the subgroup generated by $\{ga-a \mid g\in G,\,a\in
\Gamma^G_{i-1}(A)\}$.  We say that the \emph{action is nilpotent}
if, for some $i$, $\Gamma^G_i(A)=\{0\}$.

\begin{proposition}\label{prop:E* nilpotent}
For any nilpotent finite complex $X$, $\mathcal{E}^*(X)$ is a
nilpotent group.
\end{proposition}

\begin{proof} We shall prove that $\mathcal{E}^*(X)$ acts nilpotently on
$H_*(X)$.  Then it follows from \cite{D-Z} that $\mathcal{E}^*(X)$
is a nilpotent group.

Let $f\in\mathcal{E}^*(X)$ and consider the diagram
with exact rows obtained from the universal coefficient
theorem
\begin{displaymath}
\xymatrix{0 \ar[r] & \mathrm{Ext}(H_{*-1}(X),\mathbb{Z}) \ar[r]
\ar[d]^{\mathrm{Ext}(f_{*},\mathbb{Z})} & H^*(X) \ar[r]
\ar[d]^{f^*=1} &\mathrm{Hom}(H_*(X),\mathbb{Z})
\ar[r]\ar[d]^{\mathrm{Hom}(f_*,\mathbb{Z})}& 0 \\ 0 \ar[r] &
\mathrm{Ext}(H_{*-1}(X),\mathbb{Z})\ar[r] & H^*(X) \ar[r]
&\mathrm{Hom}(H_*(X),\mathbb{Z})\ar[r]& 0.\\}
\end{displaymath}
Then both $\mathrm{Ext}(f_{*},\mathbb{Z})$ and
$\mathrm{Hom}(f_*,\mathbb{Z})$ are identity maps. Write
$H_*(X)=F\oplus T$ as the sum of its free and torsion parts and
let $p_T \colon F\oplus T\to T$ and $p_F \colon F\oplus T\to F$ be
the projections.  Since
$\mathrm{Ext}(\mathbb{Z}/m,\mathbb{Z})=\mathbb{Z}/m$, it follows
that $\mathrm{Ext}(H_*(X),\mathbb{Z})=
\mathrm{Ext}(T,\mathbb{Z})=T$. Thus
$\mathrm{Ext}(f_{*},\mathbb{Z})=1 $ implies that $p_T\circ
{f_*}|_T\colon T\to T$ is the identity.

In the same way, since $\mathrm{Hom}(f_*,\mathbb{Z})=1$, we have
that $p_F\circ {f_*}|_F\colon F\to F$ is also the identity.
Therefore, $f_*\colon F\oplus T\to F\oplus T$ can be written as
$f_*(x,y)=(x,y+\phi(x))$, with $\phi\colon F\to T$ a homomorphism
that depends on $f$. Hence $\Gamma^{\mathcal{E}^*(X)}_1(H_*(X))$
is generated by elements of the form $f_*(x,y)-(x,y)=(0,\phi(x))$.
On these elements, any $g\in\mathcal{E}^*(X)$ satisfies
$g_*(0,\phi(x))=(0,\phi(x))$. Hence
$\Gamma^{\mathcal{E}^*(X)}_2(H_*(X))=\{0\}$, and the result
follows.
\end{proof}

\providecommand{\bysame}{\leavevmode\hbox
to3em{\hrulefill}\thinspace}

\end{document}